\documentclass[article,12pt]{article}
\usepackage[a4paper, total={6in, 9in}]{geometry}
\usepackage{graphicx}
\usepackage{amsmath}
\usepackage{amssymb}
\usepackage{amsfonts}
\usepackage{latexsym}
\usepackage{url}
\usepackage{color}
\usepackage{multicol}
\usepackage{booktabs}
\usepackage[normalem]{ulem}
\usepackage{rotating}

\newenvironment{keywords}{\begin{quote}\small \em 
{\bf Keywords\/}:}{\end{quote}}

\usepackage{authblk}

\begin{document}
%
\title{A comparison of formulations for aircraft deconfliction}
%
%
\author{Renan Spencer Trindade\thanks{This publication was supported by the Chair ``Integrated Urban Mobility'', backed by L’X - \'Ecole Polytechnique and La Fondation de l’\'Ecole Polytechnique. The Partners of the Chair shall not under any circumstances accept any liability for the content of this publication, for which the author shall be solely liable.} }
\author{Claudia D'Ambrosio}
%
%
\affil{LIX, CNRS, École Polytechnique, Institut Polytechnique de Paris, 
Route de Saclay, Palaiseau, 91128, France \\
\normalfont{\texttt{
\{rst,dambrosio\}@lix.polytechnique.fr}
}}

\date{}
\renewcommand\Affilfont{\itshape\small}
\maketitle              
\begin{abstract}
In this work, we aim to compare different methods and formulations to solve a problem in air traffic management to global optimality. In particular, we focus on the aircraft deconfliction problem, where we are given $n$ aircraft, their position at time $0$, and their (straight) trajectories. We wish to identify and solve potential pairwise conflict by temporarily modifying the aircraft’s trajectory. A pair of aircraft is in conflict when they do not respect a minimum, predefined safety distance. In general, conflicts could be solved both varying the aircraft's speed or trajectory, but in this paper we only consider the latter, more precisely heading-angle deviations.
The problem has been formulated as a mixed integer nonlinear program (MINLP). We compare this formulation, solved by open-source MINLP solvers for global optimization, against a reformulation that shows a larger number of variables and constraints but only separable nonconvexities. We solve such a separable formulation with the same MINLP solvers or the Sequential Convex Mixed Integer Nonlinear Programming method. The separable formulation, despite being larger, facilitates some solvers in finding good-quality solutions.
\end{abstract}

\section{Introduction}

The identification and solution of aircraft conflicts during en-route flights is a critical issue in air traffic management, since maintaining a separation between each pairs of aircraft sharing the same air sector is essential for air traffic safety. The growing volume of air traffic is increasing the challenges for air traffic controllers in guaranteeing safety and emphasizes the need for decision-making tools helping them in this process, technically called aircraft deconfliction. In the literature, considerable research attention has been devoted to developing mathematical models and efficient algorithms for aircraft deconfliction.

\begin{figure}[h]
\centering
\includegraphics[width=.55\textwidth]{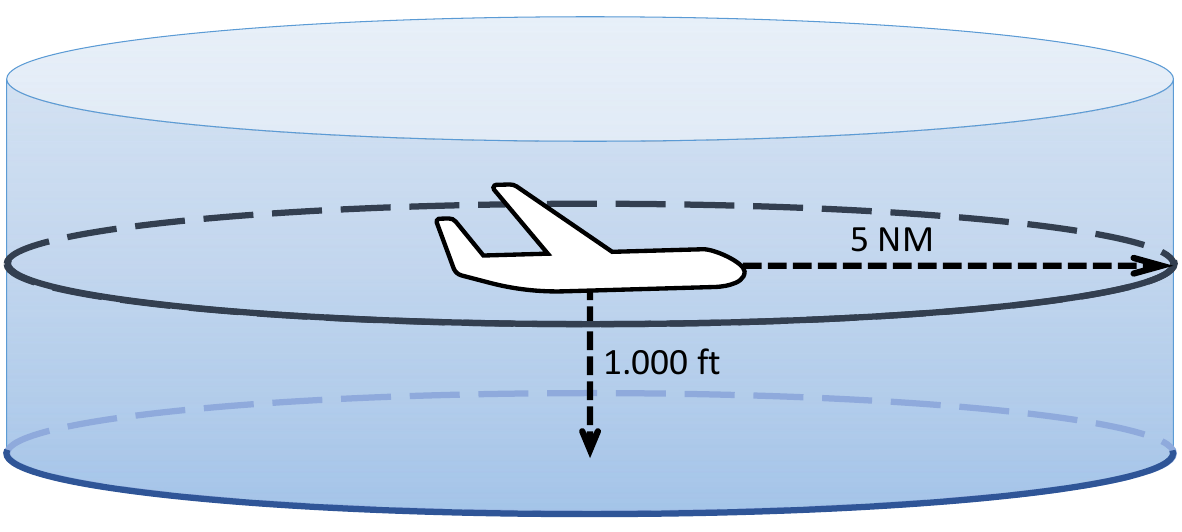}
\caption{Safety distance cylinder.} \label{fig1}
\end{figure}

Conflicts between aircraft arise when they come too close to each other, violating the required safety distances based on their predicted trajectories, see Figure \ref{fig1}. Detecting and solving these conflicts is crucial for ensuring air traffic safety. Air traffic controllers typically rely on maneuvers that involve changing the aircraft's velocity, trajectory (heading angle), or flight level (altitude) to achieve separation.

In this paper, we focus on models that only consider heading angle deviation, which is currently the most used by air traffic controllers.
The goal is to minimize deviations from the original aircraft flight plan, thus minimizing the impact on flight efficiency while guaranteeing the safety distance is satisfied by each pair of aircraft at any time.

Several papers present formulations for this problem using Mixed Integer Linear Programs (MILPs), which include time discretization \cite{Richards2002}, disjunctive linear separation constraint \cite{Pallottino2002}, and space discretization in \cite{Omer2015} for 2D conflict resolution problems. Velocity and altitude changes are considered in \cite{Alonso-Ayuso2011}. Minimizing the fuel expended is explored considering speed control and flight-level assignment considering potential conflict points \cite{Vela2009}, and linear approximations \cite{Vela2010}.
Mixed Integer Non Linear Programs (MINLPs) appear naturally well-suited for addressing the problem. Furthermore, the advances on the performances of MINLP solver now makes it possible to think of using them in real-world applications as the aircraft deconfliction. An overview of these approaches can be found in \cite{Cafieri2017_1}, where the formulation strategies allow the simultaneous consideration of continuous and discrete decision-making variables.
The work \cite{Cafieri2017} proposes a two-step solution approach that solves first a velocity control problem and then a heading angle change in sequence.
A feasibility pump heuristic is found in \cite{Cafieri2018}.
In \cite{Rey2018} is presented a complex number formulation that considers both speed and heading control. 
In \cite{Cerulli2019,DBLP:journals/jgo/CerulliDLP21}, the problem is solved with bilevel programming using formulations based on KKT conditions.
An extensive review of mathematical representations of aircraft separation conditions and presented in a unifying analysis can be found in \cite{DBLP:journals/transci/PelegrinD22}.

This paper starts presenting the MINLP formulation published in \cite{Cafieri2017}, which shows continuous and binary variables, a quadratic objective function, and nonlinear constraints. Then, in Section \ref{sec:3}, we introduce a reformulation of this model where all nonlinear constraints are expressed as univariate constraints. In Section \ref{sec:4}, we compare the formulations using Couenne \cite{Belotti2009}, SCIP \cite{BestuzhevaEtal2021OO} and SCMINLP \cite{DLW09,DAmbrosio2012,DAmbrosio2019}. We end the paper with some conclusions and perspectives in Section \ref{sec:5}.

\section{Aircraft Deconfliction Problem}
\label{sec:2}

In the aircraft deconfliction problem, we consider $A$, a set of $n$ aircraft, flying at the same altitude and sharing the same air sector within the same time window. The pairwise conflicts, i.e., situation in which pairs of aircraft are too close to one another, should be ``deconflicted'', namely their trajectories have to be modified so that the conflict is solved.
In this section, we present a formulation introduced in \cite{Cafieri2017}, using concepts first presented in \cite{Cafieri2014}. It is assumed that the initial position of each aircraft $i \in A$ is represented by the coordinates $(x_{i}^0,y_{i}^0)$, as well as the velocity $v_i$, and the starting heading angle $\phi_i$. All the mentioned data are known \emph{a priori} and provided (no uncertainty is considered). 

The decision variables $\theta_i$ represent the heading angle deviation for each aircraft $i \in A$, and it is bounded by the interval $[\theta_i^{min}, \theta_i^{max}]$, assuming the uniform motion laws applied to each aircraft.
Let us define $B$ the set of pairs of aircraft, i.e., $B := \{(i,j) \in A \times A:i<j\}$. Given a pair of aircraft $(i,j) \in B$, their relative initial distance $X_{ij}^{r_0}$ and their relative speed $V_{ij}^r$ are the following vectors in $\mathbb{R}^2$:

$$X_{ij}^{r_0} := \left(\begin{array}{@{}c@{}}
    x_{i}^0-x_{j}^0 \\
    y_{i}^0-y_{j}^0 \\
    \end{array} \right)
$$

$$
V_{ij}^r := \left(\begin{array}{@{}c@{}}
    \cos(\phi_i + \theta_i)v_i -\cos(\phi_j + \theta_j)v_j \\
    \sin(\phi_i + \theta_i)v_i -\sin(\phi_j + \theta_j)v_j  \\
    \end{array} \right).
$$

Let us introduce as well the variables $t_{ij}^m$ (for $(i,j) \in B$), which define the time when the relative distance for the pair of aircraft $i$ and $j$ is minimal. If this value is negative, aircraft are diverging, 
thus no conflict is possible within the considered time horizon 
(unless they are in conflict at instant 0, but, in this case, the instance is infeasible, thus we can remove it from our test bed).

The separation condition, i.e., the constraint that imposes that each pair of aircraft $(i,j) \in B$ are distant at least $d$ any time within the considered time horizon, can be written as follow (see \cite{Cafieri2017} for details on how to derive such constraints):
$$y_{ij}(||V_{ij}^r||^2(||X_{ij}^{r_0}||^2-d^2)-(X_{ij}^{r_0} \cdot V_{ij}^r)^2) \geq 0$$
where $X_{ij}^{r_0} \cdot V_{ij}^r$ is the inner product of $X_{ij}^{r_0}$ and $V_{ij}^r$ and the binary variable $y_{ij}$ is introduced to deactivate the the separation condition when not needed, namely when $t^m_{ij}$ is negative. Thus, we also have the following constraints
$$t_{ij}^m(2y_{ij}-1) \geq 0,$$
which imposes that $y_{ij}$ is 1 when $t_{ij}^m$ is positive and that $y_{ij}$ is 0 if $t_{ij}^m$ is negative.

The complete formulation for the aircraft deconfliction problem, introduced in \cite{Cafieri2017}, is defined below:
\begin{align} \hspace{0cm}
  (\mbox{M}_1)\;&\min_{\theta, t^m, y, V^r} \; \sum_{i \in A} \theta^2_i & \label{eq:m1:fo}\\
  \mbox{s.t.:} & \nonumber \\ 
    & \theta_i^{min} \leq \theta_i \leq \theta_i^{max} & \forall i \in A \label{eq:m1:bound}\\
    & y_{ij}(||V_{ij}^r||^2(||X_{ij}^{r_0}||^2-d^2)-(X_{ij}^{r_0} \cdot V_{ij}^r)^2) \geq 0 & \forall (i,j) \in B \label{eq:m1:distance}\\
    & t_{ij}^m = - \frac{X_{ij}^{r_0} \cdot V_{ij}^r}{||V_{ij}^r||^2} & \forall (i,j) \in B \label{eq:m1:time}\\
    & t_{ij}^m(2y_{ij}-1) \geq 0, &  \forall (i,j) \in B  \label{eq:m1:activation}\\
    & V_{ij}^r := \left(\begin{array}{@{}c@{}}
    \cos(\phi_i + \theta_i)v_i -\cos(\phi_j + \theta_j)v_j \\
    \sin(\phi_i + \theta_i)v_i -\sin(\phi_j + \theta_j)v_j  \\
    \end{array} \right)                         & \forall (i,j) \in B \label{eq:m1:V_def} \\
    & y_{ij} \in \{0,1\}  						& \forall (i,j) \in B.
  \end{align}
The objective function \eqref{eq:m1:fo} minimizes the heading angle variations for all aircraft, which is bounded by \eqref{eq:m1:bound} for each aircraft. 
The constraints \eqref{eq:m1:distance} ensure the minimum separation distance between each pair of aircraft. This condition is activated only when $t_{ij}^m \geq 0$, as guaranteed by constraints \eqref{eq:m1:time} and \eqref{eq:m1:activation}. Finally, constraints \eqref{eq:m1:V_def} model the definition of variables $V_{ij}^r$ as functions of variables $\theta$.

Note that (M$_1$) is a MINLP as it shows continuous and binary variables, quadratic objective function, and nonlinear constraints.

\section{Univariate reformulation}
\label{sec:3}
This paper aims to compare the performance of open-source MINLP solvers. Some of these use commercial ones to solve subproblems, however, they are all available for academic use. We consider the formulation proposed in \cite{Cafieri2017}. As stated there, there are a few open-source solvers which can deal with trigonometric functions. The authors used only Couenne \cite{Belotti2009} because at the time of the publication of \cite{Cafieri2017}, SCIP could not deal with trigonometric functions. However, since then, this functionality of SCIP was added, specifically starting from version 8.0 \cite{BestuzhevaEtal2021OO}. Moreover, the solver SCMINLP (see \cite{DLW09,DAmbrosio2012,DAmbrosio2019}) was recently made available at the following link: \url{https://github.com/iumx-chair/SC-MINLP}. It deals with MINLPs with non-convexities that are sums of univariate functions. 
%
Therefore, for including SCMINLP in the comparison, we derive a reformulation of model (M$_1$), where all nonlinear constraints are univariate. 

Note that, we consider only open-source solvers for MINLP which are non-convex, aiming at finding a global optimum or a valid bound to it.

\subsection{Expanding and separating the constraints}

We start the reformulation with constraints \eqref{eq:m1:distance}, which present the minimum distance separation condition between pairs of aircraft $(i,j) \in B$.
%
Let us introduce the following parameters 
to simplify the notation:
\begin{align}
&    C_{ij}=||X_{ij}^{r_0}||^2-d^2 & \forall (i,j) \in B \label{eq:Cij_def}\\
&    D_{ij}=x_{i}^0-x_{j}^0 & \forall (i,j) \in B        \\
&    E_{ij}=y_{i}^0-y_{j}^0 & \forall (i,j) \in B        \\
&    H_{ij} =
     C_{ij} (v_i^2 + v_j^2) & \forall (i,j) \in B.
\end{align}

Furthermore, let us introduce two sets of variables that represent the interaction between the head angle changes of two aircraft:
\begin{align}
    & \Phi^-_{ij}= (\phi_i + \theta_i)-(\phi_j + \theta_j) & \forall (i,j) \in B\\
    & \Phi^+_{ij}= (\phi_i + \theta_i)+(\phi_j + \theta_j) & \forall (i,j) \in B.
\end{align}

Note that, all the new notation introduced in this section is according to the following rules: capital letters stay for parameters, and capital Greek letters for variables. For the notation introduced in the previous session, we decided to use the one presented in \cite{Cafieri2017} for coherence, thus it will not follow the rules mentioned here.

Let us now consider the left-hand-side (LHS) of \eqref{eq:m1:distance}, but temporarily neglecting the multiplication per $y_{ij}$. After expanding it, we can separate and isolate the terms with only one variable used, using trigonometric properties (see the appendix \ref{appendix} for a detailed derivation). In this way, we can obtain the following univariate constraints:

\begin{align}
\Gamma_{ij} \leq & 
    - (\cos(\phi_i + \theta_i)D_{ij} v_i)^2 
    - (\sin(\phi_i + \theta_i)E_{ij} v_i)^2 \nonumber \\ 
    & - \cos(\phi_i + \theta_i) \sin(\phi_i + \theta_i)2 D_{ij} E_{ij} v_i^2 & \forall (i,j) \in B
\label{eq:42}
\end{align}

\begin{align}
    \Delta_{ij} \leq &
    - (\cos(\phi_j + \theta_j)D_{ij} v_j )^2 
    - (\sin(\phi_j + \theta_j)E_{ij} v_j )^2 \nonumber \\
    & - \cos(\phi_j + \theta_j) \sin(\phi_j + \theta_j)2D_{ij}E_{ij} v_j^2 &
    \forall (i,j) \in B
\label{eq:43}
\end{align}

\begin{align}
    & \Lambda^-_{ij} \leq
    \cos(\Phi^-_{ij}) v_i v_j ({D_{ij}}^2 + {E_{ij}}^2 - 2 C_{ij}) & \forall (i,j) \in B \\
    & \Lambda^+_{ij} \leq 
    \cos(\Phi^+_{ij})v_i v_j ({D_{ij}}^2 - {E_{ij}}^2) 
    + \sin(\Phi^+_{ij})2v_i v_j D_{ij} E_{ij} & \forall (i,j) \in B \label{eq:44}
\end{align}

Thus, we can rewrite the LHS, excluding the multiplication per $y_{ij}$, as a sum of univariate constraints:
\begin{align}
    & \Gamma_{ij} + \Delta_{ij} + \Lambda^-_{ij} + \Lambda^+_{ij} + H_{ij}  & \forall (i,j) \in B \label{eq:distance}
\end{align}


\subsection{Activation constraints}

The activation constraints \eqref{eq:m1:activation} define the variable $y_{ij}=1$ only when $t_{ij}^m$ is positive, and 0 otherwise. First, we replace $t^m_{ij}$ with its definition \eqref{eq:m1:time}, obtaining:
$$- \frac{X_{ij}^{r_0} \cdot V_{ij}^r}{||V_{ij}^r||^2} (2y_{ij}-1) \geq 0.$$
Since the value of $||V_{ij}^r||^2$ is always positive, we can ignore this division, which does not influence the sign of the LHS. Thus, we are left with:
$$- X_{ij}^{r_0} \cdot V_{ij}^r (2y_{ij}-1) \geq 0.$$
To separate this constraint into univariate constraints, we can rewrite $- X_{ij}^{r_0} \cdot V_{ij}^r$ as the sum of two additional variables $\Omega^-_{ij} + \Omega^+_{ij}$:

\begin{align}
    & \Omega^-_{ij} = -D_{ij}\cos(\phi_i + \theta_i)v_i - E_{ij}\sin(\phi_i + \theta_i)v_i & \forall (i,j) \in B\\
    & \Omega^+_{ij} = D_{ij}\cos(\phi_j + \theta_j)v_j + E_{ij}\sin(\phi_j + \theta_j)v_j & \forall (i,j) \in B.
\end{align}
The binary variable $y$ will be reintroduced in the next section within a BigM constraint involving $\Omega^-_{ij} + \Omega^+_{ij}$.

\subsection{BigM strategy}

To avoid the activation of the separation constraint by multiplying the binary variable $y_{ij}$ and keep our formulation univariate, we adopt the BigM strategy. Therefore, the activation constraints \eqref{eq:m1:activation} can be rewritten as follows:

\begin{align}
    & -M^-_{ij} (1-y_{ij}) \leq \Omega^-_{ij} + \Omega^+_{ij} \leq M^+_{ij} y_{ij} & \forall (i,j) \in B.
\end{align}

The minimum separation condition constraint can also be rewritten using the same strategy as follows:

\begin{align}
    & \Gamma_{ij} + \Delta_{ij} + \Lambda^-_{ij} + \Lambda^+_{ij} + H_{ij} \geq M_{ij} (1-y_{ij}) & \forall (i,j) \in B.
\end{align}

To define all the $M$ values, we calculate bounds for each univariate function \eqref{eq:42}--\eqref{eq:44}. To do this, an option would be to set lower and upper bound values to -1 and +1, respectively, for each sine and cosine functions and compute the smallest (or largest) value accordingly. However, we get better values for $M$ by computing the minimum values (subject to simple bounds on the variables) of each univariate function as follows:

\begin{align*}
    & M_{ij} = \min( \Gamma_{ij} ) + \min( \Delta_{ij} ) + \min( \Lambda^-_{ij} ) + \min( \Lambda^+_{ij}) + H_{ij} & \forall (i,j) \in B \\
    & M^-_{ij} = \min(\Omega^-_{ij}) + \min(\Omega^+_{ij}) & \forall (i,j) \in B \\
    & M^+_{ij} = \max(\Omega^-_{ij}) + \max(\Omega^+_{ij}) & \forall (i,j) \in B.
\end{align*}

Note that, even if these computations look expensive, MINLP solvers often solve costly optimization problems in the presolving phase, see, for example, optimization-based bound tightening. The average CPU time needed to find each of the values of bigMs is 0.007 (maximum 0.4).

\subsection{Final formulation}
Finally, we transformed the objective function into a linear one by introducing a new variable $\Theta$, and linking it to the original objective function thanks to a new inequality. The complete formulation can be found below:

\begin{align} \hspace{0cm}
  (\mbox{M}_2)\;&\min \; \Theta & \label{Model1_FO}\\
  \mbox{s.t.:} & \nonumber \\ 
    & \Theta \geq \sum_{i \in A} \theta^2_i\\
    & \Phi^-_{ij}= (\phi_i + \theta_i)-(\phi_j + \theta_j) & \forall (i,j) \in B\\
    & \Phi^+_{ij}= (\phi_i + \theta_i)+(\phi_j + \theta_j) & \forall (i,j) \in B\\
    & \Lambda^-_{ij} \geq
    \cos(\Phi^-_{ij}) v_i v_j ({D^2_{ij}} + {E^2_{ij}} - 2 C_{ij}) & \forall (i,j) \in B \\
    & \Lambda^+_{ij} \geq 
    \cos(\Phi^+_{ij})v_i v_j ({D^2_{ij}} - {E^2_{ij}}) 
    + \sin(\Phi^+_{ij})2v_i v_j D_{ij}E_{ij} & \forall (i,j) \in B \\ 
    & \Gamma_{ij} + \Delta_{ij} + \Lambda^-_{ij} + \Lambda^+_{ij} + H_{ij} \geq M_{ij}(1-y_{ij}) & \forall (i,j) \in B \label{eq:distance_bis}\\
    & -M^-_{ij}(1-y_{ij}) \leq \Omega^-_{ij} + \Omega^+_{ij} \leq M^+_{ij} y_{ij} & \forall (i,j) \in B \\
    & \Omega^-_{ij} = -D_{ij}\cos(\phi_i + \theta_i)v_i - E_{ij}\sin(\phi_i + \theta_i)v_i & \forall (i,j) \in B\\
    & \Omega^+_{ij} = D_{ij}\cos(\phi_j + \theta_j)v_j + E_{ij}\sin(\phi_j + \theta_j)v_j & \forall (i,j) \in B \\
    & \theta_i^{min} \leq \theta_i \leq \theta_i^{max} & \forall i \in A \\
    & y_{ij} \in \{0,1\}  						& \forall (i,j) \in B.
  \end{align}

\section{Computational results}
\label{sec:4}

\begin{figure}[ht]
\centering
\begin{minipage}{.40\textwidth}
\includegraphics[width=\textwidth]{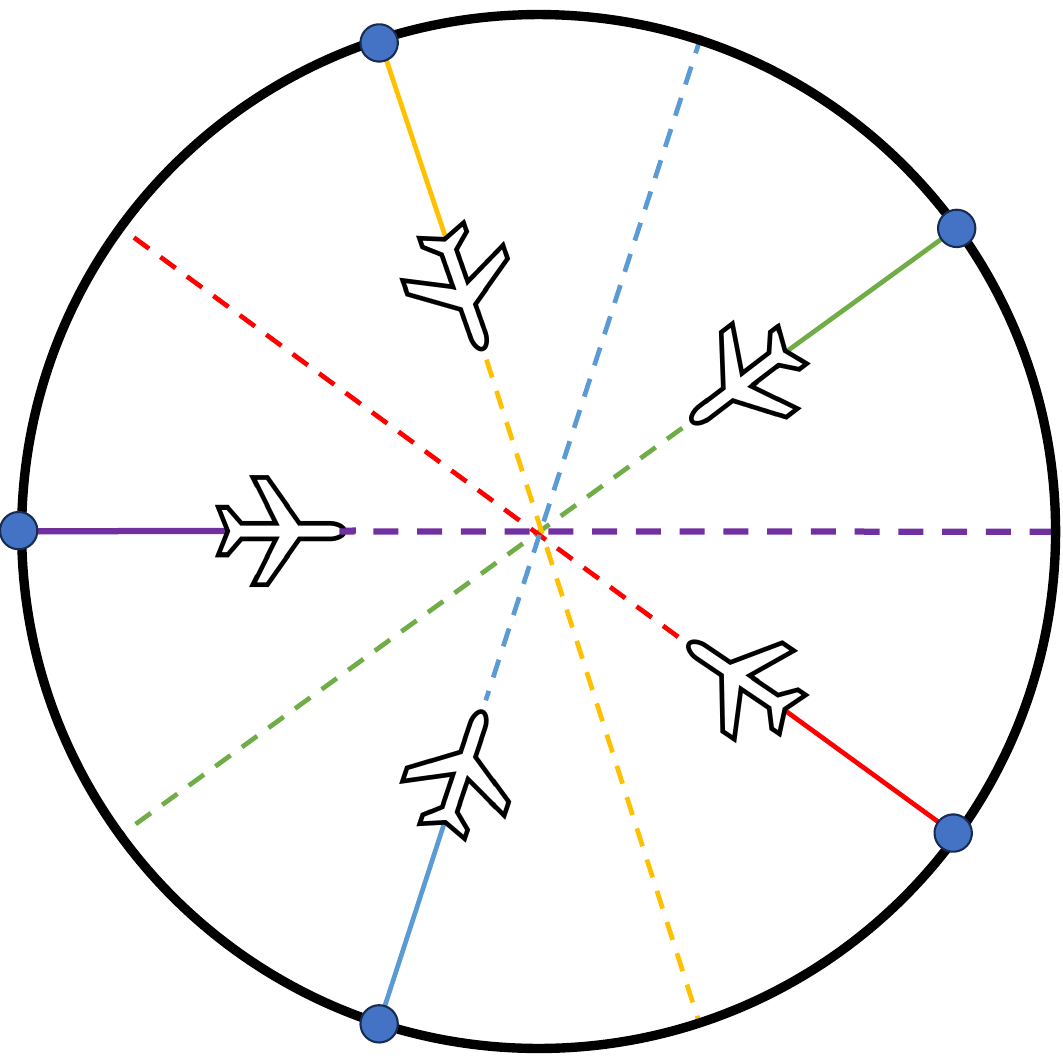}
\caption{Example of CP instances.}
\label{fig:ex2_1}
\end{minipage}%
\hspace{.01\textwidth}
\begin{minipage}{.40\textwidth}
\includegraphics[width=\textwidth]{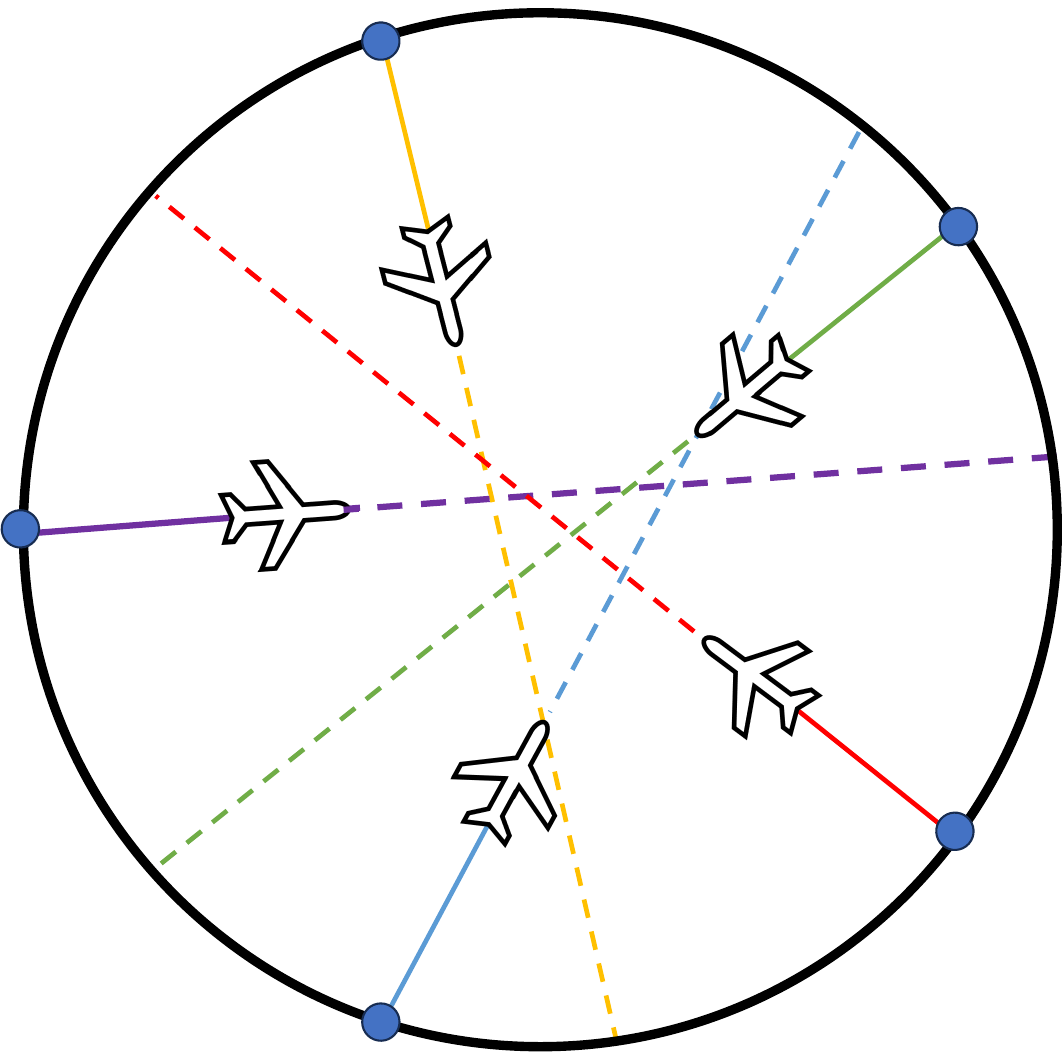}
\caption{Example of RCP instances.}
\label{fig:ex2_2}
\end{minipage}
\end{figure}

We perform our computational tests on the public instances reported in \cite{Rey2018}, which comprises two sets of instances. The first is the Circle Problem (CP), which contains a set of aircraft positioned uniformly around a circle, with the head angles oriented to the center (see, for example, Figure \ref{fig:ex2_1}). The second, called Random Circle Problem (RCP), is generated similarly, but the initial velocities and directions are randomly shifted within specified intervals  (see, for example, Figure \ref{fig:ex2_2}).

The two models were implemented in the AMPL language and tested with three solvers: COUENNE \cite{Belotti2009}, SCIP 8.0.3 \cite{BestuzhevaEtal2021OO}, and SCMINLP \cite{DLW09,DAmbrosio2012,DAmbrosio2019}, on an Intel(R) Xeon(R) CPU E5-2620 v4 @ 2.10GHz computer, with 64GB RAM.
All solvers use the default configuration with a time limit of 600 seconds.
Both SCIP and SCMINLP were run using the external libraries of CPLEX 12.10 \cite{cplex1210}.

Table \ref{res:tab1} shows the computational results for the CP instances. For each of the 18 instances (identified by the number of aircraft, varying from 3 to 20), we report the results for:
\begin{itemize}
\item Couenne\_Ori: Couenne on the original model;
\item Couenne\_Ref: Couenne on the reformulated model;
\item SCIP\_Ori: SCIP on the original model;
\item SCIP\_Ref: SCIP on the reformulated model;
\item SCMINLP: the solver SCMINLP on the reformulated model.
\end{itemize}
For each solver/model combination, we report the total processing time (Time), the objective function value of the best feasible solution found within the time limit, if any (Primal), and the dual bound (Dual) at the end of the solver execution.

\begin{table}[ht!] \hspace{-15px}
\setlength{\tabcolsep}{3.5pt}
\scalebox{0.7}{
\begin{tabular}{r|rrr|rrr|rrr|rrr|rrr}
\hline
   & \multicolumn{3}{c|}{Couenne\_Ori}              & \multicolumn{3}{c|}{Couenne\_Ref}              & \multicolumn{3}{c|}{SCIP\_Ori}                       & \multicolumn{3}{c|}{SCIP\_Ref}                      & \multicolumn{3}{c}{SCMINLP}                     \\
 $|A|$  & Time          & Primal & Dual  & Time          & Primal & Dual  & Time    & Primal & Dual  & Time    & Primal & Dual  & Time    & Primal & Dual  \\
 \hline
 3  & 0.664         & 0.001  & 0.001 & 1.257         & 0.001  & 0.001 & 1.210   & 0.001  & 0.001 & 0.750   & 0.001  & 0.001 & 67.260  & 0.001  & 0.001 \\
 4  & 3.799         & 0.001  & 0.001 & 3.064         & 0.001  & 0.001 & 600.000 & 0.001  & 0.001 & 1.840   & 0.001  & 0.001 & 277.230 & 0.001  & 0.001 \\
 5  & 99.033        & 0.002  & 0.002 & 600.865       & 0.002  & 0.001 & 600.000 & 0.002  & 0.002 & 29.040  & 0.002  & 0.002 & 600.000 & 0.002  & 0.002 \\
 6  & 602.320       & 0.004  & 0.001 & 601.301       & 0.004  & 0.000 & 41.620  & 0.004  & 0.004 & 600.000 & 0.004  & 0.004 & 600.010 & 0.004  & 0.000 \\
 7  & 603.069       & 0.006  & 0.000 & 601.449       & 0.006  & 0.000 & 600.000 & 0.006  & 0.006 & 600.000 & 0.006  & 0.002 & 600.010 & 0.006  & 0.001 \\
 8  & 602.800       & 0.035  & 0.000 & 601.845       & 0.011  & 0.000 & 600.000 & 0.011  & 0.006 & 600.000 & 0.011  & 0.001 & 600.010 & 0.011  & 0.000 \\
 9  & 603.307       & 0.020  & 0.000 & 601.656       & 0.012  & 0.000 & 600.000 & 0.012  & 0.005 & 600.000 & 0.012  & 0.000 & 600.020 & 0.012  & 0.000 \\
 10 & 602.680       & 0.139  & 0.000 & 601.342       & 0.017  & 0.000 & 600.000 & 0.017  & 0.003 & 600.000 & 0.017  & 0.000 & 600.020 & 0.017  & 0.000 \\
 11 & 601.529       & 0.260  & 0.000 & 601.532       & 0.025  & 0.000 & 600.000 & 0.049  & 0.000 & 600.000 & 0.022  & 0.000 & 600.010 & 0.022  & 0.000 \\
 12 & 601.298       & -      & 0.000 & 601.236       & 0.639  & 0.000 & 600.010 & 1.587  & 0.000 & 600.000 & 0.028  & 0.000 & 600.040 & 0.028  & 0.000 \\
 13 & 600.559       & -      & 0.000 & 599.958       & 0.664  & 0.000 & 600.000 & 0.037  & 0.000 & 600.000 & 0.037  & 0.000 & 600.050 & 0.037  & 0.000 \\
 14 & 600.264       & -      & 0.000 & 599.716       & 0.664  & 0.000 & 600.000 & 1.648  & 0.000 & 600.000 & 0.044  & 0.000 & 600.040 & 0.044  & 0.000 \\
 15 & 600.072       & -      & 0.000 & 599.123       & 0.527  & 0.000 & 600.000 & 1.873  & 0.000 & 600.000 & 0.059  & 0.000 & 600.080 & 0.059  & 0.000 \\
 16 & 600.099       & -      & 0.000 & 594.433       & 0.824  & 0.000 & 600.000 & 0.181  & 0.000 & 600.000 & 0.085  & 0.000 & 600.130 & 0.085  & 0.000 \\
 17 & 599.749       & -      & 0.000 & 597.860       & 0.538  & 0.000 & 600.000 & 0.179  & 0.000 & 600.000 & 0.093  & 0.000 & 600.070 & 0.093  & 0.000 \\
 18 & 599.651       & -      & 0.000 & 596.131       & 0.780  & 0.000 & 600.000 & 3.727  & 0.000 & 600.000 & 0.097  & 0.000 & 600.180 & 0.097  & 0.000 \\
 19 & 599.242       & -      & 0.000 & 595.040       & 1.180  & 0.000 & 600.000 & 0.233  & 0.000 & 600.000 & 0.123  & 0.000 & 600.160 & 0.124  & 0.000 \\
 20 & 599.030       & -      & 0.000 & 589.511       & 1.415  & 0.000 & 600.000 & 4.213  & 0.000 & 600.000 & 0.132  & 0.000 & 600.160 & 0.132  & 0.000 \\
 \hline
\end{tabular}
}
\caption{Results for the Circle Problem (CP) instances.}
\label{res:tab1}
\end{table}

We can see that only small instances (up to 6 aircraft) can be solved to optimality by at least one of these approaches. 
SCIP with the original model shows the best results concerning the dual bound, outperforming the other approaches on 4 instances ($|A| \in \{7,10\}$). For smaller instances, other methods could find the same bound. However, for instances with more than 10 aircraft, all the approaches found the trivial bound 0.

From what concerns the primal bound, Couenne cannot even find any feasible solution for instances with more at least 12 aircraft, when provided with the original formulation. However, it could find a feasible solution for all the instances when provided with the reformulated model.

In general, the approaches with the reformulated model show significantly better results when we analyze the quality of the primal solutions. We can see that although the dual bounds are worse, the primal solutions are equal or better compared to the original model. The approach finding the best feasible solutions is SCIP on the reformulated model. However, it is worth mentioning that SCMINLP found exactly the same solution of SCIP besides for 19 aircraft (0.124 instead of 0.123).

Table \ref{res:tab2} shows the computational results for the RCP instances. In this case, the number of aircraft varies in the set $\{10, 20, 30, 40\}$. For each of these values, 100 instances are available, thus we present aggregated results. For each pair solver/model, we report: the sum of the dual bounds (Dual bound), the number of solutions solved to optimality (Solv.), and the number of instances with feasible solutions (Feas.).

\begin{table}[ht!] \hspace{-1cm}
\setlength{\tabcolsep}{3.5pt}
\scalebox{0.7}{
\begin{tabular}{r|rrr|rrr|rrr|rrr|rrr}
\hline
& \multicolumn{3}{c|}{Couenne\_Ori} & \multicolumn{3}{c|}{Couenne\_Ref}              & \multicolumn{3}{c|}{SCIP\_Ori} & \multicolumn{3}{c|}{SCIP\_Ref}               & \multicolumn{3}{c}{SCMINLP}  \\
$|A|$ & Dual bound        & Solv. & Feas.  & Dual bound         & Solv.   & Feas.  & Dual bound         & Solv. & Feas.  & Dual bound         & Solv. & Feas.  & Dual bound         & Solv. & Feas.  \\
\hline
10       & 5,007E-04          & 9         & \textbf{100} & 4,454E-05          & \textbf{14} & \textbf{100} & \textbf{5,847E-02} & 9         & \textbf{100} & 3,693E-02          & 6         & \textbf{100} & 1,815E-03          & 1         & \textbf{100} \\
20       & -4,866E-08         & 0         & 0            & 0,000E+00          & 0           & 35           & 0,000E+00          & 0         & 99           & 1,637E-07          & 0         & \textbf{100} & \textbf{1,340E-03} & 0         & 76           \\
30       & 0,000E+00          & 0         & 0            & 0,000E+00          & 0           & 0            & 0,000E+00          & 0         & \textbf{56}  & 0,000E+00          & 0         & 19           & \textbf{3,589E-04} & 0         & 3            \\
40       & \textbf{0,000E+00} & 0         & 0            & \textbf{0,000E+00} & 0           & 0            & \textbf{0,000E+00} & 0         & \textbf{10}  & \textbf{0,000E+00} & 0         & 0            & -1,053E-04         & 0         & 0            \\
\hline
\end{tabular}
}
\caption{Results for the Random Circle Problem (RCP) instances.}
\label{res:tab2}
\end{table}


We could note that, for 10 aircraft, all the methods could find a feasible solution for each of the instances. However, Couenne applied to the reformulated model shows largest number of instances solved to global optimality. Starting from 20 aircraft instances, none of the approaches could find a global optimum and prove global optimality. However, it is clear that SCIP outperfoms the other methods when it comes to finding feasible solutions (in particular, when run on the original model). Note also that SCMINLP performs better than Couenne on both models.

For what concerns the dual bound, SCIP on the original model and SCMINLP perform better than the other approaches.

The next tables are introduced to perform a comparison on the the quality of the feasible solutions.
In particular, in Table \ref{res:tab4} we report in each cell the number of RCP instances (over 400) for which both the approach on the row and the column could find a feasible solution. This information is relevant because we next compare pairwise the methods only on these instances.

\begin{table}[]
\centering
\setlength{\tabcolsep}{3.5pt}
\begin{tabular}{l|c|c|c|c|c}
\hline
             & Couenne\_Ori & Couenne\_Ref & SCIP\_Ori & SCIP\_Ref & SCMINLP \\
\hline
Couenne\_Ori & - & 100                      & 100                      & 100                      & 100                      \\
Couenne\_Ref & 100                      & - & 135                      & 135                      & 127                      \\
SCIP\_Ori    & 100                      & 135                      & - & 211                      & 178                      \\
SCIP\_Ref    & 100                      & 135                      & 211                      & - & 176                      \\
SCMINLP      & 100                      & 127                      & 178                      & 176                      & - \\
\hline
\end{tabular}
\caption{Number of RCP instances for which both approaches find a feasible solution.}
\label{res:tab4}
\end{table}

In Table \ref{res:tab3}, we report in each cell the sum of the objective function value of the feasible solution found by the method in the row on the instances for which both this method and the one in the column could find a feasible solution.

\begin{table}[ht!]
\centering
\setlength{\tabcolsep}{3.5pt}
\begin{tabular}{l|c|c|c|c|c}
\hline
             & Couenne\_Ori & Couenne\_Ref & SCIP\_Ori & SCIP\_Ref & SCMINLP \\
\hline
Couenne\_Ori & - & \textbf{0.215}                    & 3.127                    & 1.731                    & 1.731                    \\
Couenne\_Ref & 1.731                    & - & 28.173                   & 28.173                   & 22.473                   \\
SCIP\_Ori    & \textbf{2.026}                    & 66.708                   & - & 202.642                  & 135.323                  \\
SCIP\_Ref    & \textbf{0.151}                    & \textbf{1.982}                   & \textbf{39.773}                   & - & 3.847                    \\
SCMINLP      & \textbf{0.377}                    & \textbf{1.360}                    & \textbf{3.578}                    & \textbf{3.414 }                 & - \\
\hline
\end{tabular}
\caption{Comparing the quality of the feasible solutions found by both approaches on the RCP instances.}
\label{res:tab3}
\end{table}

For example, in the cell identified by the first row and second column, we report the sum of the objective function value of the feasible solutions found by Couenne on the original model, considering only the 100 instances solved by Couenne on both the original and reformulated model. This value is 0.215, which is smaller than 1.731, namely the sum of the objective function value of the feasible solutions found by Couenne on the reformulated model on the same set of instances. Having highlighted in bold the best performances, we can clearly see that SCMINLP consistently shows the best performance from what concerns solution quality, on the instances for which it could find a feasible solution. Let us remind, however, the reader that SCMINLP could find less feasible solutions than SCIP. Couenne is clearly dominated on both aspects by both SCMINLP and SCIP. The last observation is important as, in the literature, the great majority of the mathematical optimization approaches relying on the use of Couenne -- as earlier versions of SCIP could not deal with trigonometric functions.

Finally, we report in Table \eqref{res:tab5} the number of instances of the RCP test bed for which the approach on the row found a feasible solution and the one on the column did not, showing that there is not a full dominance of SCIP with respect to SCMINLP as the latter could find a feasible solution for 3 (resp. 1) instances for which SCIP on the reformulated (resp. original) model could not find any. The same can be said for Couenne on the reformulated model, as it could find feasible solutions for 8 instances for which SCMINLP could not find any. However, Couenne on the original model is dominated by all the other approaches.

\begin{table}[]
\centering
\setlength{\tabcolsep}{3.5pt}
\begin{tabular}{l|c|c|c|c|c}
\hline
             & Couenne\_Ori & Couenne\_Ref & SCIP\_Ori & SCIP\_Ref & SCMINLP \\
\hline
Couenne\_Ori & -            & 0            & 0         & 0         & 0       \\
Couenne\_Ref & 35           & -            & 0         & 0         & 8       \\
SCIP\_Ori    & 165          & 130          & -         & 54        & 87      \\
SCIP\_Ref    & 119          & 84           & 8         & -         & 43      \\
SCMINLP      & 79           & 52           & 1         & 3         & -      \\
\hline
\end{tabular}
\caption{Number of RCP instances for which the approach on the row found a feasible solution and the one on the column did not.}
\label{res:tab5}
\end{table}

\section*{Conclusions}
\label{sec:5}
In this paper, we computationally analyzed two MINLP formulations for the aircraft deconfliction problem, with heading-angle deviations. We did not consider speed changes to focus on the impact of sin/cos operators. In particular, we were interested in testing the performance of SCIP, which was recently generalized to deal with such operators, and SCMINLP, a solver for MINLP problems with separable non-convesities, recently distributed open-source.

The first considered formulation is taken from \cite{Cafieri2017}. The second is a separable reformulation, i.e., a formulation where each non linearity is represented by a sum of univariate functions. The latter was introduced to include SCMINLP to our computational comparison.

We considered instances from the literature 418 instances from the literature, with a number of aircraft varying from 3 to 40. On these, the solvers that provides the best dual bounds are either SCIP on the formulation \cite{Cafieri2017} or SCMINLP (on the separable formulation), the latter in particular for instances with 20 and 30 aircraft. For the largest instances (40 aircraft)all solvers are stuck with the trivial bound of $0$.

SCIP on the formulation \cite{Cafieri2017} is also the most performing concerning the number of instances for which a feasible solution was found (283). However, in a pairwise comparison considering only the instances for which both solvers found a solution, SCMINLP outperforms the others, and SCIP on the separable reformulation outperforms all the approaches besides SCMINLP.

To conclude, as SCIP is now dealing with sin/cos operators, it looks very interesting to use it for solving the aircraft deconfliction problem, for which up to now Couenne was used in the literature. The newly distributed solver SCMINLP could be interesting for finding good quality solutions. Methods combining the advantages of the two solvers could be interesting to explore.

%
%
%
\bibliographystyle{splncs04}
\bibliography{numta}

\appendix
\section*{Appendix}

\section{Deriving \eqref{eq:distance}}\label{appendix}
We start from
\begin{equation}
    ||V_{ij}^r||^2(||X_{ij}^{r_0}||^2-d^2))-(X_{ij}^{r_0}\cdot  V_{ij}^r)^2
\end{equation}
which can be written as
\begin{equation}
    C_{ij} ||V_{ij}^r||^2-(X_{ij}^{r_0} \cdot V_{ij}^r)^2 \label{eq:appendix1}
\end{equation}
thanks to \eqref{eq:Cij_def}.
%
%

We now start with the expansion of the first term of \eqref{eq:appendix1}:
\begin{multline}
    C_{ij}||V_{ij}^r||^2=\\ C_{ij}(\cos(\phi_i + \theta_i)v_i -\cos(\phi_j + \theta_j)v_j)^2 + \\
    C_{ij}(\sin(\phi_i + \theta_i)v_i -\sin(\phi_j + \theta_j)v_j)^2 = \\ 
%
%
%
    C_{ij}[(\cos(\phi_i + \theta_i)v_i)^2 + \\
    - \cos(\phi_i + \theta_i) \cos(\phi_j + \theta_j)2v_i v_j + \\
    + (\cos(\phi_j + \theta_j)v_j)^2 + \\
    + (\sin(\phi_i + \theta_i)v_i)^2 + \\
    - \sin(\phi_i + \theta_i) \sin(\phi_j + \theta_j)2v_j v_i + \\
    + (\sin(\phi_j + \theta_j)v_j)^2] =  \\
%
    - \cos(\phi_i + \theta_i) \cos(\phi_j + \theta_j)2 C_{ij} v_i v_j + \\
    - \sin(\phi_i + \theta_i) \sin(\phi_j + \theta_j)2 C_{ij} v_j v_i + \\
    + C_{ij} v_i^2 + C_{ij} v_j^2 \\
\label{eq:132}
\end{multline}

Now we can consider the second term of \eqref{eq:appendix1}
\begin{multline}
    (X_{ij}^{r_0}\cdot  V_{ij}^r)^2=\\
    [D_{ij}(\cos(\phi_i + \theta_i)v_i -\cos(\phi_j + \theta_j)v_j)+ \\
    (E_{ij}(\sin(\phi_i + \theta_i)v_i -\sin(\phi_j + \theta_j)v_j )]^2 \hspace{2.65cm} = \\
%
%
%
%
%
    [\cos(\phi_i + \theta_i)D_{ij}v_i -\cos(\phi_j + \theta_j)D_{ij}v_j)+ \\
    (\sin(\phi_i + \theta_i)E_{ij}v_i -\sin(\phi_j + \theta_j)E_{ij}v_j )]^2 = \\
%
%
%
    + (\cos(\phi_i + \theta_i)D_{ij}v_i)^2 + \\
    - \cos(\phi_i + \theta_i) \cos(\phi_j + \theta_j)2{D^2_{ij}} v_iv_j + \\
    + (\cos(\phi_j + \theta_j)D_{ij}v_j)^2 + \\
    + \cos(\phi_i + \theta_i) \sin(\phi_i + \theta_i)2D_{ij}E_{ij}v_i^2 + \\
    - \cos(\phi_j + \theta_j) \sin(\phi_i + \theta_i)2D_{ij}E_{ij} v_i v_j + \\
    + (\sin(\phi_i + \theta_i)E_{ij}v_i)^2 + \\
    - \cos(\phi_i + \theta_i) \sin(\phi_j + \theta_j)2D_{ij}E_{ij} v_jv_i + \\
    + \cos(\phi_j + \theta_j) \sin(\phi_j + \theta_j)2D_{ij}E_{ij} v_j^2 + \\
    - \sin(\phi_i + \theta_i) \sin(\phi_j + \theta_j)2 {E^2_{ij}} v_j v_i + \\
    + (\sin(\phi_j + \theta_j)E_{ij}v_j)^2 \\
\label{eq:22}
\end{multline}

%
%
%
From \eqref{eq:22}, we define $\Gamma_{ij}$ (resp. $\Delta_{ij}$) the sum of terms involving exclusively the variable $\theta_i$ (resp. $\theta_j$), see \eqref{eq:42} (resp. \eqref{eq:43}).
%
%
%
%
%
%

Now, subtracting from \eqref{eq:132} what is left of \eqref{eq:22}, we have the following multivariate expression to reformulate:
\begin{multline}
- \cos(\phi_i + \theta_i) \cos(\phi_j + \theta_j)2 v_i v_j (C_{ij} - {D^2_{ij}}) \\
- \sin(\phi_i + \theta_i) \sin(\phi_j + \theta_j)2 v_i v_j (C_{ij} - {E^2_{ij}}) + \\
+ \cos(\phi_i + \theta_i) \sin(\phi_j + \theta_j)2 v_i v_j D_{ij}E_{ij} + \\
+ \sin(\phi_i + \theta_i) \cos(\phi_j + \theta_j)2 v_i v_j D_{ij}E_{ij} + \\
+ C_{ij} v^2_i + C_{ij} v^2_j \\
\label{eq:33}
\end{multline}


Now we can use the trigonometric identities to convert into univariate functions the first 4 terms of the equation \eqref{eq:33}. First we introduce two additional set of variables (for $(i,j) \in B$):

\begin{align}
    \Phi^-_{ij}= (\phi_i + \theta_i)-(\phi_j + \theta_j) \label{eq:50}\\
    \Phi^+_{ij}= (\phi_i + \theta_i)+(\phi_j + \theta_j). \label{eq:51}
\end{align}

Then, we apply the well-known trigonometric identities:

\begin{multline}
\cos(\phi_i + \theta_i) \cos(\phi_j + \theta_j) =
\frac{1}{2} ( \cos(\Phi^-_{ij}) + \cos(\Phi^+_{ij}) )\\
\sin(\phi_i + \theta_i) \sin(\phi_j + \theta_j) =
\frac{1}{2} ( \cos(\Phi^-_{ij}) - \cos(\Phi^+_{ij}) )\\ 
\cos(\phi_i + \theta_i) \sin(\phi_j + \theta_j) =
\frac{1}{2} ( \sin(\Phi^+_{ij}) - \sin(\Phi^-_{ij}) )\\ 
\sin(\phi_i + \theta_i) \cos(\phi_j + \theta_j) =
\frac{1}{2} ( \sin(\Phi^-_{ij}) + \sin(\Phi^+_{ij}) )\\ 
\label{eq:52}
\end{multline}


Now we can simplify and separate into monovariate functions:

\begin{align}
\Lambda^-_{ij} = & \cos(\Phi^-_{ij}) v_i v_j ({D_{ij}}^2 + {E_{ij}}^2 - 2 C_{ij} )  \\ 
\Lambda^+_{ij} = & \cos(\Phi^+_{ij})v_i v_j ({D_{ij}}^2 - {E_{ij}}^2) + \sin(\Phi^+_{ij})2v_i v_j D_{ij}E_{ij}  \\ 
H_{ij} = & C_{ij} v_i^2 + C_{ij} v_j^2 
\label{eq:52ter}
\end{align}

The final function is given by the sum of the univariated functions.

\begin{equation}
    C_{ij}||V_{ij}^r||^2 - (X_{ij}^{r_0}\cdot  V_{ij}^r)^2=
    \Gamma_{ij} + \Delta_{ij} + \Lambda^-_{ij} + \Lambda^+_{ij} + H_{ij}
\label{eq:55}
\end{equation}



\end{document}